\documentclass[12pt]{article}

\usepackage{amsmath,amssymb,amsthm,amscd,latexsym}
\usepackage[colorlinks, linkcolor=blue, anchorcolor=gree, citecolor=red]{hyperref}
\usepackage[numbers,sort&compress]{natbib}
\usepackage{graphicx}

\usepackage{mathtools}
\usepackage{color}
\usepackage[normalem]{ulem}

\newtheorem{thm}{Theorem}[section]

\newtheorem{lem}[thm]{Lemma}

\newtheorem{cor}[thm]{Corollary}
\theoremstyle{definition}
\newtheorem{ex}[thm]{Example}

\newtheorem{remark}[thm]{Remark}
\def\thefootnote{\fnsymbol{footnote}}
\newcommand{\bth}{\begin{thm}}
\renewcommand{\eth}{\end{thm}}
\newcommand{\bex}{\begin{examp}}
\newcommand{\eex}{\end{examp}}
\newcommand{\bre}{\begin{remark}}
\newcommand{\ere}{\end{remark}}

\newcommand{\bal}{\begin{aligned}}
\newcommand{\eal}{\end{aligned}}
\newcommand{\beq}{\begin{equation}}
\newcommand{\eeq}{\end{equation}}
\newcommand{\ben}{\begin{equation*}}
\newcommand{\een}{\end{equation*}}

\newcommand{\bpf}{\begin{proof}}
\newcommand{\epf}{\end{proof}}
\renewcommand{\thefootnote}{}

\renewcommand{\thefootnote}{\arabic{footnote}}
\def\beql#1{\begin{equation}\label{#1}}

\definecolor{VeryLightBlue}{rgb}{0.9,0.9,1}
\definecolor{LightBlue}{rgb}{0.8,0.8,1}
\definecolor{MidBlue}{rgb}{0.5,0.5,1}
\definecolor{DarkBlue}{rgb}{0,0,0.6}
\definecolor{Blue}{rgb}{0,0,1}

\definecolor{Gold}{rgb}{1,0.843,0}

\definecolor{LightGreen}{rgb}{0.88,1,0.88}
\definecolor{MidGreen}{rgb}{0.6,1,0.6}
\definecolor{DarkGreen}{rgb}{0,0.6,0}

\definecolor{VeryLightYellow}{rgb}{1,1,0.9}
\definecolor{LightYellow}{rgb}{1,1,0.6}
\definecolor{MidYellow}{rgb}{1,1,0.5}
\definecolor{DarkYellow}{rgb}{1,1,0.2}

\definecolor{DarkPurple}{rgb}{.6,0,1}

\definecolor{Red}{rgb}{1,0,0}
\definecolor{VeryLightRed}{rgb}{1,0.9,0.9}
\definecolor{LightRed}{rgb}{1,0.8,0.8}
\definecolor{MidRed}{rgb}{1,0.55,0.55}

\def\Cay{{\rm Cay}}
\def\FFF{\mathbb{F}}
\def\ZZZ{\mathbb{Z}}

\def\Ga{\Gamma}

\def\qed{\hfill$\Box$\vspace{11pt}}

\long\def\delete#1{}

\begin{document}

\title{\Large\bf Subgroup perfect codes in Cayley graphs}

\author{\renewcommand{\thefootnote}{\arabic{footnote}}Xuanlong Ma\footnotemark[1] , Gary L. Walls\footnotemark[2] , Kaishun Wang\footnotemark[3] , Sanming Zhou\footnotemark[4]}

\footnotetext[1]{School of Science, Xi'an Shiyou University, Xi'an 710065, China}

\footnotetext[2]{Department of Mathematics, Southeastern Louisiana University, Hammond, LA 70402, USA}

\footnotetext[3]{Laboratory of Mathematics and Complex Systems (Ministry of Education), School of Mathematical Sciences, Beijing Normal University, Beijing 100875, China}

\footnotetext[4]{School of Mathematics and Statistics, The University of Melbourne, Parkville, VIC 3010, Australia}

\renewcommand{\thefootnote}{}
\footnotetext{E-mail addresses: \texttt{xuanlma@mail.bnu.edu.cn} (X. Ma), \texttt{gary.walls@selu.edu} (G. L. Walls), \texttt{wangks@bnu.edu.cn} (K. Wang), \texttt{sanming@unimelb.edu.au} (S. Zhou)}

\date{}

\maketitle
\newcommand\blfootnote[1]{%
\begingroup
\renewcommand\thefootnote{}\footnote{#1}%
\addtocounter{footnote}{-1}%
\endgroup
}

\begin{abstract}
Let $\Gamma$ be a graph with vertex set $V(\Gamma)$.
A subset $C$ of $V(\Gamma)$ is called a perfect code in $\Gamma$ if $C$ is an independent set of $\Gamma$ and every vertex in $V(\Gamma)\setminus C$ is adjacent to exactly one vertex in $C$. A subset $C$ of a group $G$ is called a perfect code of $G$ if there exists a Cayley graph of $G$ which admits $C$ as a perfect code. A group $G$ is said to be code-perfect if every proper subgroup of $G$ is a perfect code of $G$. In this paper we prove that a group is code-perfect if and only if it has no elements of order $4$. We also prove that a proper subgroup $H$ of an abelian group $G$ is a perfect code of $G$ if and only if the Sylow $2$-subgroup of $H$ is a perfect code of the Sylow $2$-subgroup of $G$. This reduces the problem of determining when a given subgroup of an abelian group is a perfect code to the case of abelian $2$-groups. Finally, we determine all subgroup perfect codes in any generalized quaternion group.

\medskip
{\em Keywords:} Perfect code; Efficient dominating set; Subgroup perfect code; Cayley graph; Finite group

\medskip
{\em AMS subject classifications (2010):} 05C25, 05C69, 94B25
\end{abstract}

\section{Introduction}

Perfect codes are important objects of study in coding theory ever since the beginning of information theory. Roughly speaking, a code is perfect if it achieves maximum possible error correction without ambiguity. In the classical setting, much work has been focused on perfect codes under the Hamming or Lee metric. Solving a long-standing conjecture, it was proved in the 1970s \cite{Va73, Tie73, ZL73} that the well-known Hamming and Golay codes are the only nontrivial linear perfect codes under the Hamming metric. (A \emph{linear code} is a subspace of some linear space $\FFF_q^n$, where $\FFF_q$ is the field with $q$ elements, $q$ being a prime power and $n$ a positive integer.) In contrast to the linear case, there are many nonlinear perfect codes under the Hamming metric, and the study of them has long been an active research area in coding theory. The reader is referred to the survey papers \cite{Heden1, Va75} for a large number of results on perfect codes under the Hamming metric. With regard to the Lee metric, the famous Golomb-Welch conjecture asserts that for any $n > 2$, $e > 1$ and $q \ge 2t+1$  there is no $q$-ary perfect $t$-codes of length $n$ under the Lee metric. A central problem for Lee codes, this 50-year-old conjecture is still wide open \cite{HK18} despite extensive research on the topic.

From a mathematical point of view, perfect codes can be defined for any finite metric space: Given an integer $t \ge 1$, a subset of a finite metric space is called a \emph{perfect $t$-code} \cite{Va75} if the balls of radius $t$ with centres in the subset form a partition of the space. In particular, since any graph is a metric space under the usual graph distance, we can talk about perfect $t$-codes in graphs. Let $\Gamma$ be a graph with vertex set $V(\Gamma)$ and edge set $E(\Gamma)$. The \emph{distance} in $\Ga$ between two vertices $u, v \in V(\Gamma)$, denoted by $d(u, v)$, is the length of a shortest path between $u$ and $v$ in $\Gamma$, and is defined to be $\infty$ if no path between $u$ and $v$ exists. In view of the definition above, a subset $C$ of $V(\Gamma)$ is a \emph{perfect $t$-code} \cite{Kr86} in $\Gamma$ if every vertex of $\Gamma$ is at distance no more than $t$ to exactly one vertex of $C$. In what follows a perfect $1$-code is simply called a {\em perfect code}. It is readily seen that a subset $C$ of $V(\Gamma)$ is a perfect code in $\Gamma$ if and only if $C$ is an independent set of $\Gamma$ and every vertex in $V(\Gamma) \setminus C$ is adjacent to exactly one vertex in $C$. In graph theory, a perfect code in a graph is also called an efficient dominating set \cite{DeS} or independent perfect dominating set \cite{Le} of the graph.

The {\em Cartesian product} of $n$ graphs $\Ga_1, \Ga_{2}, \ldots, \Ga_{n}$ is the graph with vertex set $V(\Ga_1) \times V(\Ga_{2}) \times \cdots \times V(\Ga_{n})$ such that two vertices $(u_1, u_2, \ldots, u_n), (v_1, v_2, \ldots, v_n)$ are adjacent if and only if $u_{i} \neq v_{i}$ for exactly one subscript $i$, and for this $i$, $u_{i}$ and $v_{i}$ are adjacent in $\Ga_i$. The \emph{Hamming graph} $H(n, q)$ is the Cartesian product of $n$ copies of the complete graph $K_q$ with $q$ vertices. Denote by $C_q^{\Box n}$ the Cartesian product of $n$ copies of the cycle $C_q$ of length $q$. In particular, $H(n, 2)$ is the $n$-dimensional cube $Q_n$ and $C_q^{\Box 2}$ is the grid graph on a torus. Alternatively, we may define $H(n, q)$ and $C_q^{\Box n}$ to be the graphs with vertex set $\ZZZ_q^n$ such that two elements $(u_1, u_2, \ldots, u_n), (v_1, v_2, \ldots, v_n)$ of $\ZZZ_q^n$ are adjacent in $H(n, q)$ if and only if they differ at exactly one coordinate and adjacent in $C_q^{\Box n}$ if and only if $u_{i} \neq v_{i}$ for exactly one $i$ and moreover $u_{i} \equiv v_{i} \pm 1 \mod q$ for this $i$.

It is well known that the Hamming and Lee metrics over $\ZZZ_q^n$ are exactly the graph distances in $H(n, q)$ and $C_q^{\Box n}$, respectively. Therefore, perfect $t$-codes under the Hamming or Lee metric are exactly those in $H(n, q)$ or $C_q^{\Box n}$, respectively. It is also well known that all Hamming graphs are distance-transitive. (A graph $\Ga$ is called \emph{distance-transitive} if for any $u, v, u', v' \in V(\Ga)$ with $d(u, v) = d(u', v')$ there exists an automorphism of $\Ga$ which maps $(u, v)$ to $(u', v')$.) This motivated Biggs \cite{Big} to study perfect codes in distance-transitive graphs as a generalization of perfect codes under the Hamming metric. Among other things he generalized the celebrated Lloyd's Theorem \cite{Lens} to perfect codes in any distance-transitive graph. The seminal paper of Biggs \cite{Big} and the fundamental work of Delsarte \cite{Del} inspired much work on perfect codes in distance-transitive graphs and, in general, in distance-regular graphs and association schemes. It is known \cite{Chi, MZ95} that many infinite families of classical distance-regular graphs have no nontrivial perfect codes, including the Grassmann graphs and the bilinear forms graphs. Doob graphs are an important family of distance-regular Cayley graphs. A necessary and sufficient condition for a Doob graph to admit perfect codes was recently given in \cite{Kro20}, and all possible parameters of subgroup perfect codes in Doob graphs were described in \cite{SHK19}. It is known \cite{Best} that if $q$ is not a prime power then a perfect $t$-code in $H(n,q)$ exists only when $t \in \{1, 2, 6, 8\}$, but the existence of such codes is a long-standing open problem. More results on perfect codes in distance-regular graphs can be found in, for example, \cite{HS, SE, Thas80}.

\medskip
\textbf{Perfect codes in Cayley graphs.}~
As observed in \cite{HXZ18}, perfect codes in Cayley graphs are another generalization of perfect codes in the classical setting. This is so because $H(n, q)$ and $C_q^{\Box n}$ are both Cayley graphs of $\ZZZ_q^n$. In general, given a group $G$ with identity element $e$ and an inverse-closed subset $S$ of $G$ with $e \notin S$, the {\em Cayley graph} ${\rm Cay}(G, S)$ of $G$ with \emph{connection set} $S$ is defined to be the graph with vertex set $G$ such that two distinct elements $x,y$ are adjacent if and only if $yx^{-1}\in S$, where a subset $A$ of $G$ is called {\em inverse-closed} if $A^{-1}:=\{a^{-1}: a \in A\} = A$. Obviously, ${\rm Cay}(G, G \setminus \{e\})$ is the complete graph with vertex set $G$, and ${\rm Cay}(G, \emptyset)$ is the graph on $G$ with no edges.

In recent years, perfect codes in Cayley graphs have received considerable attention \cite{D14, DSLW16, FHZ, HXZ18, KM13, Le, MBG07, OPR07, Ta13, Z15}. The reader is referred to \cite[Section 1]{HXZ18} for a brief account of results on perfect codes in Cayley graphs and connections between such codes and factorizations and tilings of the underlying groups. In general, a \emph{tiling} \cite{Dinitz06} of a group $G$ is a pair of subsets $(A, B)$ of $G$ such that $e \in A \cap B$ and every element of $G$ can be expressed uniquely as $ab$ with $a \in A$ and $b \in B$. It is readily seen that $(A, B)$ is a tiling of $G$ such that $A$ is inverse-closed if and only if $B$ is a perfect code of $\Cay(G, A \setminus \{e\})$ such that $e \in B$.

In \cite{HXZ18}, Huang, Xia and Zhou introduced the following concept: A subset $C$ of a group $G$ is called a {\em perfect code} of $G$ if there exists a Cayley graph $\Cay(G, S)$ of $G$ which admits $C$ as a perfect code. In particular, a perfect code of $G$ which is also a subgroup of $G$ is called a {\em subgroup perfect code} of $G$. In the same paper, Huang, Xia and Zhou obtained a necessary and sufficient condition for a normal subgroup of a group $G$ to be a perfect code of $G$, and determined all subgroup perfect codes of all dihedral groups and some abelian groups. As explained in \cite{HXZ18}, in some sense subgroup perfect codes are an analogue of linear perfect codes.

\medskip
\textbf{Code-perfect groups.}~
It may happen that every subgroup of a given group is a perfect code. We call a group with this property a code-perfect group. More explicitly, a group $G$ is said to be \emph{code-perfect} if for every subgroup $H$ of $G$ there exists a Cayley graph $\Cay(G, S)$ of $G$ which admits $H$ as a perfect code (that is, $(S \cup \{e\}, H)$ is a tiling of $G$ with $S \cup \{e\}$ inverse-closed). Note that the trivial subgroup $\{e\}$ is a perfect code in the complete Cayley graph ${\rm Cay}(G, G \setminus \{e\})$ and the whole group $G$ is a perfect code in the empty Cayley graph ${\rm Cay}(G, \emptyset)$. So a code-perfect group can also be defined as a group in which every proper subgroup is a perfect code in the group.

It is natural to ask which groups are code-perfect. In this paper, we answer this question by giving a complete characterization of code-perfect groups. As we will see shortly, by this characterization not all abelian groups are code-perfect. So one may ask when a given subgroup of an abelian group is a perfect code. We give an answer to this question by reducing the problem of determining when a subgroup of an abelian group is a perfect code to the case of abelian $2$-groups. It turns out that no generalized quaternion group can be code-perfect. We determine all subgroup perfect codes together with corresponding Cayley graphs in any generalized quaternion group.

\medskip
\textbf{Notation and terminology.}~Before stating our results let us introduce some notation and terminology first. All groups considered in the paper are finite, and all graphs considered are finite and undirected with no loops or multiple edges. So we will omit the adjective ``finite'' before the words ``group'' and ``graph''. We always use $e$ to denote the identity element of the group under consideration. We use $G_2$ and $G_{2'}$ to denote the Sylow $2$-subgroup and Hall $2'$-subgroup of a group $G$, respectively. Note that, for any abelian group $G$, $G_2$ consists of the elements of $G$ with order a power of $2$, and $G_{2'}$ consists of the elements of $G$ with odd order. Denote
$$
G^2 = \{g^2: g \in G\}.
$$
An abelian group $G$ is said to be \emph{$2$-divisible} if $G^2 = G$. For an abelian group $G$, a subgroup $H$ of $G$ is called a {\em $2$-pure subgroup} of $G$ if $G^2 \cap H=H^2$. As usual, we use $A \times B$ to denote the direct product of two groups $A$ and $B$. We use $Q_{4n}$ to denote the generalized quaternion group of order $4n$, where $n \ge 2$. It is well known (see, for example, \cite[pp. 44--45]{Jon}) that
\begin{equation}
\label{eq:q4n}
Q_{4n}=\langle x,y: x^n=y^2, x^{2n}=e, y^{-1}xy=x^{-1}\rangle
\end{equation}
and the order of $x^i y$ in $Q_{4n}$ is $4$ for $0 \le i \le n$.

\medskip
\textbf{Main results.}~The first main result in this paper is as follows.

\begin{thm}
\label{newthm1}
A group is code-perfect if and only if it has no elements of order $4$.
\end{thm}

The sufficiency of this result will be proved by construction: Given any group $G$ with no elements of order $4$ and any proper subgroup $H$ of $G$, we will construct an inverse-closed subset $S$ of $H$ with $e \notin S$ such that $\Cay(G, S)$ admits $H$ as a perfect code.

In \cite[Corollary 2.4(a)]{HXZ18}, it was proved that every normal subgroup of any group of odd order is a perfect code of the group. The following corollary of Theorem \ref{newthm1} generalizes this result from all normal subgroups to all subgroups.

\begin{cor}
\label{newthm1-cor1}
Any group of odd order is code-perfect.
\end{cor}

Theorem \ref{newthm1} also implies the following result, in which $\mathbb{Z}_2^d \times Q$ is interpreted as $Q$ when $d=0$.

\begin{cor}
\label{newthm1-cor2}
An abelian group is code-perfect if and only if it is isomorphic to $\mathbb{Z}_2^d \times Q$ for some integer $d \ge 0$ and abelian group $Q$ of odd order.
\end{cor}

All simple groups with no elements of order $4$ have been classified in \cite{Wa69}. Combining this and Theorem \ref{newthm1}, we obtain the following result.

\begin{cor}
\label{newthm1-cor3}
A simple group is code-perfect if and only if it is isomorphic to one of the following groups:
\begin{itemize}
\item[{\rm (a)}] a cyclic group of prime order;

\item[{\rm (b)}] ${\rm PSL}(2,2^e)$, $e\ge 2$;

\item[{\rm (c)}] ${\rm PSL}(2,q)$, $q\equiv \pm 3~({\rm mod}~8)$, $q\ge 5$;

\item[{\rm (d)}] a Ree group $^2G_2(3^{2n+1})$, $n\ge 1$;

\item[{\rm (e)}] the Janko group $J_1$.
\end{itemize}
\end{cor}

Theorem \ref{newthm1} implies that not every abelian group is code-perfect. So one may ask when a given subgroup of an abelian group is a perfect code of the group. The following result shows that this problem can be reduced to the case of abelian $2$-groups.

\begin{thm}
\label{abethm1}
Let $G$ be an abelian group and $H$ a proper subgroup of $G$. Then $H$ is a perfect code of $G$ if and only if $H_2$ is a $2$-pure subgroup of $G_2$, which in turn is true if and only if $H_2$ is a perfect code of $G_2$.
\end{thm}

A property which is diagonally opposite to the one of being a code-perfect group is that no nontrivial proper subgroup is a perfect code. Our next result gives all abelian non-simple groups with this property. It would be interesting if one can obtain a characterization of non-abelian groups with this property.

\begin{thm}
\label{newrem}
Let $G$ be an abelian group which is not a simple group. Then every nontrivial proper subgroup of $G$ is not a perfect code of $G$ if and only if $G$ is isomorphic to the cyclic group $\mathbb{Z}_{2^m}$ for some $m\ge 2$.
\end{thm}

In the special case when $G$ is a cyclic group, Theorems~\ref{abethm1} and \ref{newrem} together yield \cite[Corollary 2.8(a)]{HXZ18}, which asserts that a proper subgroup $H$ of a cyclic group $G$ is a perfect code of $G$ if and only if either $|H|$ or $[G:H]$ is odd.

Theorem~\ref{newthm1} also implies that $Q_{4n}$ is not code-perfect. We determine all subgroup perfect codes of $Q_{4n}$ together with corresponding Cayley graphs in the following result.

\begin{thm}
\label{newmainth2}
Let $Q_{4n}$ be the generalized quaternion group as presented in \eqref{eq:q4n}, and let $H$ be a proper subgroup of $Q_{4n}$. Then $H$ is a perfect code of $Q_{4n}$ if and only if one of the following holds:
\begin{itemize}
\item[{\rm (a)}] $H=\langle x^t\rangle$, where $t$ is a positive integer dividing $2n$ such that $\frac{2n}{t}$ is odd;

\item[{\rm (b)}] $H=\langle x^t, x^sy\rangle$, where $t \ge 3$ is an odd integer dividing $2n$ and $s$ is an integer between $0$ and $t-1$.
\end{itemize}
Moreover, if (a) occurs, then $\langle x^t\rangle$ is a perfect code of $\Cay(Q_{4n}, S)$, where
\beq
\label{eq:s1}
S = \left\{x^n, x^i, x^{-i}: 1 \le i \le (t/2)-1\right\} \cup \left\{x^i y, x^{n+i} y: 0 \le i \le (t/2)-1\right\};
\eeq
and if (b) occurs, then $\langle x^t, x^sy\rangle$ is a perfect code of $\Cay(Q_{4n}, S)$, where
\beq
\label{eq:s2}
S = \left\{x^i, x^{-i}: 1 \le i \le (t-1)/2\right\}.
\eeq
\end{thm}

\medskip
\textbf{Structure of the paper.}~We will present some preliminary results in the next section. The proofs of Theorems \ref{newthm1} and \ref{newmainth2} will be given in Sections \ref{Sec1} and \ref{Se}, respectively, and the proofs of Theorems \ref{abethm1} and \ref{newrem} will be given in Section \ref{abeg}. An example to illustrate Theorem \ref{newmainth2} will be given in Section \ref{Se}.

\section{Preliminaries}
\label{sec:pre}

We will use the following result in our proofs of Theorems \ref{abethm1} and \ref{newmainth2}.

\begin{thm}{\rm (\cite[Theorem 2.2(a)]{HXZ18})}
\label{basic}
Let $G$ be a group and $H$ a normal subgroup of $G$. Then $H$ is a perfect code of $G$ if and only if the following holds: for any $g\in G$, $g^2\in H$ implies $(gh)^2=e$ for some $h\in H$.
\end{thm}

The following lemma is an extension of \cite[Lemma 2.1]{HXZ18}, where the equivalence between the first two statements was established.

\begin{lem}
\label{mainlem0}
Let $G$ be a group and $H$ a subgroup of $G$. Let $S$ be an inverse-closed subset of $G$ such that $e \not\in S$.
The following statements are equivalent:
\begin{itemize}
\item[{\rm (a)}] $H$ is a perfect code of ${\rm Cay}(G, S)$;

\item[{\rm (b)}] $S\cup\{e\}$ is a left transversal of $H$ in $G$;

\item[{\rm (c)}] $S\cup\{e\}$ is a right  transversal of $H$ in $G$.
\end{itemize}
\end{lem}

\proof
The equivalence between (a) and (b) was proved in \cite[Lemma 2.1]{HXZ18}. It remains to prove that (b) and (c) are equivalent.

Suppose that $S\cup\{e\}$ is a left transversal of $H$ in $G$. We claim that $Ha\ne Hb$ for distinct $a, b \in S$. Suppose to the contrary that $Ha = Hb$. Then $a^{-1}H=b^{-1}H$. Since $a^{-1}, b^{-1}\in S$ (as $S$ is inverse-closed) and $S\cup\{e\}$ is a left transversal of $H$ in $G$, we deduce that $a^{-1}=b^{-1}$, but this contradicts our assumption that $a\ne b$. Therefore, $\{Hs:s\in S\}$ consists of $|S|$ right cosets of $H$ in $G$. Since $S\cup\{e\}$ is a left transversal of $H$ in $G$, we have $[G:H]=|S|+1$ and $s \notin H$ for each $s \in S$. This together with $e\notin S$ implies that $S\cup\{e\}$ is a right  transversal of $H$ in $G$. So (b) implies (c). Similarly, we can prove that (c) implies (b).
\qed

Lemma \ref{mainlem0} implies the following result.

\begin{cor}
\label{gtlem0}
Let $G$ be a group and $H$ a proper subgroup of $G$. Then $H$ is a perfect code of $G$ if and only if there exists a left or right transversal of $H$ in $G$ which contains $e$ and is inverse-closed.
In particular, if there exists an element $x \in G\setminus H$ such that $xH$ or $Hx$ is inverse-closed and contains no involutions, then $H$ is not a perfect code of $G$.
\end{cor}

Of course a left or right transversal $T$ of $H$ in $G$ contains $e$ if and only if $T \cap H = \{e\}$.

\section{Proof of Theorem \ref{newthm1}}
\label{Sec1}

We will prove two lemmas before giving the proof of Theorem \ref{newthm1}. An element $x$ of a group is called a \emph{square} if it can be expressed as $x = y^2$ for some element $y$ of the group.

\begin{lem}
\label{mainlem1}
Let $G$ be a group and $x$ an involution of $G$. Then $\langle x\rangle$ is a perfect code of $G$ if and only if $x$ is not a square of $G$.
\end{lem}

\proof
Denote $H=\langle x\rangle = \{e, x\}$. If $x$ is a square, say, $x=y^2$ for some $y\in G$, then $Hy=\{y, xy\}=\{y,y^{-1}\}$ is inverse-closed and contains no involutions. Hence, by the second statement in Corollary \ref{gtlem0}, $H$ is not a perfect code of $G$.

Now assume that $x$ is not a square. We will apply induction to construct a right transversal $T$ of $H$ in $G$ which contains $e$ and is inverse-closed. Once this is achieved, we then obtain from Corollary \ref{gtlem0} that $H$ is a perfect code of $G$.

To begin with, we process initially the coset $H$ and put $e$ into $T$ to represent $H$. Inductively, suppose that we have processed some but not all right cosets of $H$ in $G$ and selected a representative for each of them, in such a way that the set of representatives selected so far is inverse-closed. Take an element $y \in G$ which is not in any right coset already processed. (For example, when only the coset $H$ has been processed, we simply take any $y \in G \setminus H$.) According to the orders of $y$ and $xy$, we now process one, two or four right cosets of $H$ in the following way.

\smallskip
{\bf Case 1.} $y$ is an involution.
\smallskip

In this case we only process $Hy$ and put $y$ into $T$ as the representative of $Hy$. (Alternatively, if $xy$ is also an involution, we can put $xy$ but not $y$ into $T$ to represent $Hy$.)

\smallskip
{\bf Case 2.} $y$ has order greater than $2$ but $xy$ is an involution.
\smallskip

In this case we only process $Hy = Hxy$ and put $xy$ into $T$ as the representative of $Hy$. We can do so because $xy$ has not been selected, for otherwise $y$ would be in a previously processed right coset of $H$ in $G$, which is a contradiction.

Note that $Hy^{-1}=\{y^{-1}, xy^{-1}\}$ and $xy^{-1}=(yx)^{-1}$ is an involution. Hence our rules in Cases 1 and 2 applied to $Hy^{-1}$ implies that $xy^{-1}$ but not $y^{-1}$ is selected to represent $Hy^{-1}$ when processing $Hy^{-1}$ (which can take place before or after $Hy$ is processed). Therefore, the undesired situation where $y^{-1}$ is a representative but $y$ is not (or $y$ is a representative but $y^{-1}$ is not) cannot happen.

\smallskip
{\bf Case 3.} Both $y$ and $xy$ have order greater than $2$.
\smallskip

Assume that $xy=yx$ first. Then $xy^{-1}=y^{-1}x$ and so $H(xy)^{-1}=\{y^{-1}x,xy^{-1}x\}=\{y^{-1}x,y^{-1}\}=Hy^{-1}$. Since $x$ is not a square, we have $y \neq y^{-1} x$. Hence $Hy=\{y,xy\}$ and $Hy^{-1}$ are distinct cosets. We process both $Hy$ and $Hy^{-1}$, and put $y$ and $y^{-1}$ into $T$ to represent $Hy$ and $Hy^{-1}$, respectively.

Now we assume that $xy\ne yx$. We have $Hy=\{y,xy\}$, $H(xy)^{-1}=\{y^{-1}x,xy^{-1}x\}$, $H(xy^{-1}x)^{-1}=\{xyx,yx\}$ and $H(yx)^{-1}=\{xy^{-1},y^{-1}\}$, and one can easily verify that these cosets are distinct and their union is inverse-closed. We process these four cosets and put $y, xy^{-1}x, xyx$ and $y^{-1}$ into $T$ as their representives, respectively.

After the treatment above, we have processed at least one more right coset of $H$ in $G$ and obtained a larger set of representatives. By our selection of representatives and based on the hypothesis, this larger set of representatives remains to be inverse-closed. If all right cosets of $H$ have been processed, we stop and output $T$. Otherwise we repeat the procedure above. By induction we can eventually obtain a transversal $T$ of $H$ in $G$ which contains $e$ and is inverse-closed, as required.
\qed

\begin{lem}\label{mainlem2}
Suppose that $G$ is a group with no elements of order $4$. Then for every subgroup $H$ of $G$ there exists a right transversal of $H$ in $G$ which contains $e$ and is inverse-closed.
\end{lem}

\proof
Let
$$
\Lambda_1 = \{Hx: x \in G \text{ and $Hx$ contains an element of order $2$}\} \setminus \{H\}
$$
and
$$
\Lambda_2=\{Hx: x\in G \text{ and $Hx$ contains no elements of order $2$}\}\setminus \{H\}.
$$

\smallskip
{\bf Claim 1.} Let $Hx\in \Lambda_2$. Then $(hx)^2\notin H$ for any $h\in H$.
\smallskip

Suppose for a contradiction that $(hx)^2\in H$ for some $h\in H$. We claim that there exist elements $y_1$ and $y_2$ of $G$ such that $y_1^2=e$, $y_2$ is of odd order and $hx=y_1y_2 = y_2y_1$. In fact, if $hx$ is of odd order, then we can take $y_1 = e$ and $y_2 = hx$. Assume that $hx$ is of even order. Then the order of $hx$ is of the form $2t$ for some odd integer $t = 2m - 1 \ge 1$ as $G$ has no elements of order $4$. Setting $y_1 = (hx)^{-t}$ and $y_2 = (hx)^{2m}$, we have $hx = y_1y_2 = y_2y_1$, $y_1^2 = e$, and $y_2$ has odd order $t$, as required.

By the assumption $(hx)^2\in H$ and the claim above, we obtain $(hx)^2=(y_1y_2)^2=y_2^2\in H$. Since the order of $y_2$ is odd, it follows that $y_2\in \langle y_2^2\rangle\subseteq H$. Now $Hx=(Hh^{-1})y_1y_2=Hy_1y_2=(Hy_2)y_1=Hy_1$. Hence $y_1 \in Hx$. Since $x\notin H$ (as $Hx \neq H$) and $y_1^2=e$, it follows that $y_1$ has order $2$, but this contradicts the assumption that $Hx\in \Lambda_2$.

\smallskip
{\bf Claim 2.} Let $Hx\in \Lambda_2$. Then $Hx\ne Hx^{-1}$ and $Hx^{-1}\in \Lambda_2$.
\smallskip

In fact, if $Hx=Hx^{-1}$, then $x^2\in H$ and so $(x^2 x)^2 = x^6\in H$, which contradicts Claim 1. Hence $Hx\ne Hx^{-1}$. Since $x\notin H$, we have $Hx^{-1}\ne H$. We claim that $Hx^{-1}\in \Lambda_2$. Suppose otherwise. Then $Hx^{-1}\in \Lambda_1$. So there exists $h\in H$ such that $(hx^{-1})^2=e$, that is, $xh^{-1} = h x^{-1}$. We then have $(h^{-1}x)^2=h^{-1}(xh^{-1})x = h^{-1}(hx^{-1})x = e\in H$, which contradicts Claim 1. Hence $Hx^{-1}\in \Lambda_2$.

\smallskip
{\bf Claim 3.} The operation $Hx \cdot h=H(xh)$ for $Hx \in \Lambda_2$ and $h\in H$ defines an action of $H$ on $\Lambda_2$.
\smallskip

In fact, for $Hx\in \Lambda_2$ and $h \in H$, since $Hxh=h^{-1}(Hx)h$, the set of orders of the elements in $Hx$ is the same as the set of orders of the elements in $Hxh$. In particular, like $Hx$, $Hxh$ contains no elements of order $2$. Moreover, since $x\notin H$ (as $Hx \ne H$), we have $Hxh \ne H$. Thus $Hxh\in \Lambda_2$ and the operation above defines an action of $H$ on $\Lambda_2$.

By Claim 2, whenever $Hx\in \Lambda_2$, we have $Hx^{-1} \in \Lambda_2$. Denote by ${\rm orb}_H(Hx)$ and ${\rm orb}_H(Hx^{-1})$ the orbits of $Hx$ and $Hx^{-1}$ under the action of $H$ defined in Claim 3, respectively.

\smallskip
{\bf Claim 4.} Let $Hx\in \Lambda_2$. Then ${\rm orb}_H(Hx) \cap {\rm orb}_H(Hx^{-1})=\emptyset$.
\smallskip

Suppose otherwise. Then $Hx^{-1}\in {\rm orb}_H(Hx)$ and so $Hx^{-1}=Hxh$ for some
$h\in H$. Hence $xhx\in H$. We then have $(hx)^2\in H$, but this contradicts Claim 1. Therefore, ${\rm orb}_H(Hx) \cap {\rm orb}_H(Hx^{-1})=\emptyset$.

\smallskip
{\bf Claim 5.} Let $Hx\in \Lambda_2$. Then $|{\rm orb}_H(Hx)|=[H:H\cap H^x]=|{\rm orb}_H(Hx^{-1})|$.
\smallskip

In fact, the stabilizer of $Hx$ under the action of $H$ is equal to $\{h\in H: Hxh=Hx\} =\{h\in H: xhx^{-1}\in H\} = \{h\in H: h\in H^x\} = H\cap H^x$. Hence $|{\rm orb}_H(Hx)|=[H:H\cap H^x]$. Similarly, $|{\rm orb}_H(Hx^{-1})|=[H:H\cap H^{x^{-1}}]$. However, as
$H\cap H^{x^{-1}}=(H\cap H^x)^{x^{-1}}$, we have $[H:H\cap H^x]=[H:(H\cap H^x)^{x^{-1}}]=[H:H\cap H^{x^{-1}}]$. Hence $|{\rm orb}_H(Hx)|=|{\rm orb}_H(Hx^{-1})|$.

\smallskip
{\bf Claim 6.} If $Hx \in \Lambda_2$ and $Hy\in \Lambda_2\setminus ({\rm orb}_H(Hx)\cup {\rm orb}_H(Hx^{-1}))$, then $Hy^{-1} \in \Lambda_2\setminus ({\rm orb}_H(Hx)\cup {\rm orb}_H(Hx^{-1}))$.

In fact, if $Hy^{-1} \in {\rm orb}_H(Hx^{\epsilon})$ where $\epsilon=\pm 1$, then $Hy^{-1}=Hx^{\epsilon}h$ for some $h \in H$. It follows that $x^{\epsilon}hy=h_1 \in H$. Hence $hy=x^{-\epsilon}h_1$. Therefore, $Hy=Hx^{-\epsilon}h_1 \in {\rm orb}_H(Hx^{-\epsilon})$, a contradiction.

We are now ready to construct a right transversal of $H$ in $G$ which contains $e$ and is inverse-closed. First, we put the identity element $e$ into the transversal to represent coset $H$. Then, for each $Hx \in \Lambda_1$, we choose an element of $Hx$ with order $2$ and put it into the transversal. It remains to select an appropriate representative for each coset in $\Lambda_2$.

By Claims 2, 4 and 6, there exist elements $x_1, x_2, \ldots, x_m$ of $G$ such that $\Lambda_2$ is partitioned into
$$
{\rm orb}_H(Hx_1),\ {\rm orb}_H(Hx_1^{-1}),\ {\rm orb}_H(Hx_2),\ {\rm orb}_H(Hx_2^{-1}),\ \ldots,\ {\rm orb}_H(Hx_m),\ {\rm orb}_H(Hx_m^{-1}).
$$
By Claim 5, for $1 \le i \le m$, we may assume that
$$
{\rm orb}_H(Hx_i)=\{Hx_i, Hx_i h_{i1}, Hx_i h_{i2},\ldots,Hx_i h_{ik_i}\}
$$
and
$$
{\rm orb}_H(Hx_i^{-1})=\{Hx_i^{-1},Hx_i^{-1}g_{i1},
Hx_i^{-1}g_{i2},\ldots,Hx_i^{-1}g_{ik_i}\}
$$
for some $h_{i1},\ldots,h_{ik_i} \in H$ and $g_{i1},\ldots,g_{ik_i}\in H$, where $k_i = [H:H\cap H^{x_i}] - 1$. Note that
$$
{\rm orb}_H(Hx_i)=\{Hx_i, Hg_{i1}^{-1} x_i h_{i1}, Hg_{i2}^{-1}x_i h_{i2},\ldots,Hg_{ik_i}^{-1}x_i h_{ik_i}\}
$$
and
$$
{\rm orb}_H(Hx_i^{-1})=\{Hx_i^{-1},Hh_{i1}^{-1}x_i^{-1}g_{i1},
Hh_{i2}^{-1}x_i^{-1}g_{i2},\ldots,Hh_{ik_i}^{-1}x_i^{-1}g_{ik_i}\}.
$$
So we can add
$$
x_i, g_{i1}^{-1} x_i h_{i1}, g_{i2}^{-1}x_i h_{i2},\ldots g_{ik_i}^{-1}x_i h_{ik_i}
$$
and
$$
x_i^{-1}, h_{i1}^{-1}x_i^{-1}g_{i1}, h_{i2}^{-1}x_i^{-1}g_{i2},\ldots, h_{ik_i}^{-1}x_i^{-1}g_{ik_i}
$$
to the transversal to represent the cosets in ${\rm orb}_H(Hx_i)$ and ${\rm orb}_H(Hx_i^{-1})$, respectively. Note that $(g_{ij}^{-1}x_i h_{ij})^{-1} = h_{ij}^{-1}x_i^{-1}g_{ij}$ for $1 \le i \le m$ and $1 \le j \le k_i$.

So far we have chosen a representative for each right coset of $H$ in $G$ and thus constructed a right transversal of $H$ in $G$. The construction itself ensures that this transversal contains $e$ and is inverse-closed.
\qed

\noindent {\em Proof of Theorem~{\rm\ref{newthm1}}.} Let $G$ be a group. If $G$ contains an element of order $4$, say, $y$, then $x = y^2$ is an involution and by Lemma~\ref{mainlem1}, $\langle x\rangle$ is not a perfect code of $G$. Hence $G$ is not code-perfect.

Now assume that $G$ has no elements of order $4$. By Lemma~\ref{mainlem2}, for every proper subgroup $H$ of $G$, there is a right transversal $T$ of $H$ in $G$ such that $e\in T$ and $T^{-1} = T$. Hence, by Corollary~\ref{gtlem0}, $H$ is a perfect code of $G$. Since this holds for any proper subgroup of $G$, we conclude that $G$ is a code-perfect group. This proves the sufficiency.
\qed

The proof of Lemma \ref{mainlem2} gives an algorithm for constructing a Cayley graph which admits a given proper subgroup of a group with no elements of order $4$ as a perfect code. In fact, if $S \cup \{e\}$ is the right transversal of $H$ in $G$ constructed in the proof of Lemma \ref{mainlem2}, then this Cayley graph is $\Cay(G, S)$.

\section{Proofs of Theorems \ref{abethm1} and \ref{newrem}}
\label{abeg}

The next lemma is obtained by applying Theorem \ref{basic} to abelian groups.

\begin{lem}
\label{new-2-group}
Let $G$ be an abelian group and $H$ a subgroup of $G$. Then $H$ is a perfect code of $G$ if and only if $H$ is a $2$-pure subgroup of $G$.
\end{lem}

\proof
We first prove the sufficiency. Suppose that $H$ is a $2$-pure subgroup of $G$. Then $G^2\cap H=H^2$. Thus, for any $g \in G$ with $g^2\in H$, we have $g^2\in H^2$ and hence there exists $h\in H$ such that $g^2=h^2$, which implies that $(gh^{-1})^2=e$. Now by Theorem~\ref{basic}, $H$ is a perfect code of $G$.

We next prove the necessity. Suppose that $H$ is a perfect code of $G$. Let $g\in G$ be such that $g^2\in G^2\cap H$. Since $H$ is normal in $G$ and $g^2 \in H$, by Theorem~\ref{basic} there exists $h\in H$ such that $(gh)^2=e$. So $g^2=(h^{-1})^2\in H^2$. It follows that $G^2\cap H\subseteq H^2$. Also, it  is clear that
$H^2 \subseteq G^2\cap H$. Hence $G^2\cap H=H^2$ and $H$ is a $2$-pure subgroup of $G$.
\qed

The following lemma is well known. We include its proof as we are unable to find a specific reference for it.

\begin{lem}
\label{2-pure}
Let $G$ be an abelian group and $H$ a subgroup of $G$. Then $H$ is a $2$-pure subgroup of $G$ if and only if $H_2$ is a $2$-pure subgroup of $G_2$.
\end{lem}

\proof
Suppose that $H_2$ is a $2$-pure subgroup of $G_2$. Since $G_2$ is a $2$-pure subgroup of $G$ and the property of being a $2$-pure subgroup is transitive, we obtain that $H_2$ is a $2$-pure subgroup of $G$. Since $H_{2'}$ has odd order, it is $2$-divisible and hence a $2$-pure subgroup of $G$. It follows that $H$ is a $2$-pure subgroup of $G$.

Now suppose that $H$ is a $2$-pure subgroup of $G$. We aim to prove that $H_2$ is a $2$-pure subgroup of $G_2$. Let $g \in G_2$ be such that $g^2 \in H_2$. Since $H$ is a 2-pure subgroup of $G$, there exists $h \in H$ such that $g^2 = h^2$. We can write $h= xy$ for some $x \in H_2$ and $y \in H_{2'}$. So $g^2 = (xy)^2 = x^2 y^2$ as $G$ is abelian. Hence $y^2 \in H_2 \cap H_{2'} = \{e\}$. Now that $y^2 = e$, we have $g^2 = x^2 \in H_2^2$. It follows that $G_2^2 \cap H_2 \subseteq H_2^2$ and therefore $H_2$ is a $2$-pure subgroup of $G_2$.
\qed

We are now ready to prove Theorem~\ref{abethm1}.
\medskip

\noindent {\em Proof of Theorem~{\rm\ref{abethm1}}.}
The first statement in Theorem~\ref{abethm1} follows from Lemmas \ref{new-2-group} and \ref{2-pure} immediately. Applying Lemma \ref{new-2-group} to the subgroup $H_2$ of $G_2$, we then obtain the second statement in Theorem~\ref{abethm1}.
\qed

A {\em complement} of a subgroup $H$ in a group $G$ is a subgroup $K$ of $G$ such that $G=HK$ and $H\cap K=\{e\}$. We will use the following lemma in our proof of Theorem \ref{newrem}.

\begin{lem}\label{comp}
Let $G$ be a group and $H$ a subgroup of $G$. If $H$ has a complement in $G$, then $H$ is a perfect code of $G$.
\end{lem}
\proof
Let $K$ be a complement of $H$ in $G$. It is easy to see that $K$ is a right transversal of $H$ in $G$. Of course $K$ contains $e$ and is inverse-closed. So, by Corollary~\ref{gtlem0}, $H$ is a perfect code of $G$.
\qed

\medskip

\noindent {\em Proof of Theorem~{\rm\ref{newrem}}.}
Suppose that $G\cong \mathbb{Z}_{2^m}$, where $m\ge 2$. Let $G=\langle g\rangle$ and let $H=\langle g^t\rangle$ be any nontrivial proper subgroup of $G$. It is clear that $t$ is even. It follows that $g^t\in G^2\cap H$. Note that $2|H^2|=|H|$ and $H^2$ is a subgroup of $G$. We have $g^t\notin H^2$ and hence $G^2\cap H\ne H^2$. This means that $H$ is not a $2$-pure subgroup of $G$, and so $H$ is not a perfect code of $G$ by Theorem~\ref{abethm1}.

Conversely, suppose that for any nontrivial proper subgroup $H$ of $G$, $H$ is not a perfect code of $G$. Then by Corollary~\ref{newthm1-cor1}, $G$ has even order.
If $G=H\times K$ with $|H|\ge 2$ and $|K|\ge 2$, then from Lemma~\ref{comp} it follows that $H$ is a perfect code of $G$, a contradiction. Therefore, $G$ is a cyclic group with order a power of $2$, as desired.
\qed

\section{Proof of Theorem \ref{newmainth2}}
\label{Se}

It is known that the subgroups of $Q_{4n}$ (where $n \ge 2$) are $\langle x^t\rangle$ and $\langle x^t, x^sy\rangle$, where $t$ is a positive integer dividing $2n$ and $s$ is an integer with $0\le s \le t-1$. Clearly, $\langle x^t, x^sy\rangle$ is either the cyclic group $\langle y\rangle$ or a generalized quaternion group.

We observe that, for any odd integer $q \ge 3$ and any odd integer $i$ between $1$ and $q-2$, the number $j = \frac{q-i}{2}$ is an integer between $1$ and $\frac{q-1}{2}$ satisfying $i+2j=q$.

\begin{lem}\label{g1}
Let $n \ge 2$ be an integer and $t$ a positive integer dividing $2n$. Then the proper subgroup $\langle x^t\rangle$ of $Q_{4n}$ is a perfect code of $Q_{4n}$ if and only if $\frac{2n}{t}$ is odd. Moreover, if $\frac{2n}{t}$ is odd, then $\Cay(Q_{4n}, S)$ with $S$ as given in \eqref{eq:s1} admits  $\langle x^t\rangle$ as a perfect code.
\end{lem}

\proof
Denote $G=Q_{4n}$ and $H=\langle x^t\rangle$.

Suppose that $H$ is a perfect code of $G$. Suppose to the contrary that $\frac{2n}{t}$ is even. That is, $|H|$ is even and so $x^n\in H$.
Take $g\in G\setminus\langle x\rangle$. Then $g^2=x^n\in H$.
Since $H$ is normal in $G$, by Theorem~\ref{basic} we have $(gh)^2=e$ for some $h\in H$.
Since $g\notin H$, we deduce that $gh$ is an involution.
It follows that $gh=x^n$, which implies that $g=x^nh^{-1}\in H$, a contradiction. Hence $\frac{2n}{t}$ must be odd.

Suppose that $\frac{2n}{t}=q$ is odd. Then $t$ is even and $|H|=q$. Using Theorem~\ref{basic}, we are going to prove that $H$ is a perfect code of $G$. Clearly, if $q=1$, then $H=\langle e\rangle$ is a perfect code of $G$. Assume that $q\ge 3$ in the sequel. Consider any $g\in G$ such that $g^2\in H$. Say, $g^2=x^{it}$ for some $1\le i \le q$. The order of $g$ is not equal to $4$. So we have $g\in \langle x\rangle$ and $g=x^k$ for some $1\le k \le 2n$. Hence $x^{2k}=x^{it}$. If $g\in H$, then $(gh)^2=e$ for $h = g^{-1}\in H$. Assume that $g\in \langle x\rangle\setminus H$ in the remaining proof. Set $k' = k$ if $1 \le k \le n$ and $k' = k-n$ if $n < k \le 2n$. Since $x^{2n}=e$ and $x^{2k}=x^{it}$, we have $x^{2k'}=x^{it}$. Since $1 \le k' \le n$, it follows that $2k' = it$. So $i$ is odd as $g\notin H$. If $i=q$, then $g=x^n$ and taking $h=e$ we obtain that $(gh)^2=e$. Now assume $1 \le i \le q-2$. Setting $j = \frac{q-i}{2}$, by the observation before Lemma \ref{g1} we see that $j$ is an integer satisfying $i+2j=q$ and $1 \le j \le \frac{q-1}{2}$. Since $2k' = it$ and $tq=2n$, we have $k' + j t = n$. Thus $(gh)^2=e$ for $h=x^{jt}\in H$. In summary, we have proved that for any $g\in G$ with $g^2\in H$ there exists $h \in H$ such that $(gh)^2=e$. Therefore, by Theorem~\ref{basic}, $H$ is a perfect code of $G$.

We now construct an inverse-closed subset $S$ of $G \setminus \{e\}$ such that $H=\langle x^t\rangle$ is a perfect code in $\Cay(G, S)$ under the condition that $\frac{2n}{t}=q$ is odd. If $q = 1$, then we can take $S = G \setminus \{e\}$ (which agrees with \eqref{eq:s1} as $t=2n$). Assume that $q \ge 3$ in the sequel. It is clear that $\{x^i: 0 \le i \le t-1\}$ is a right transversal of $H$ in $\langle x \rangle$.
Since $x^{2n-i} = x^{(q-1)t}x^{t-i} \in Hx^{t-i}$ for each $i$, it follows that $\left\{e, x^i, x^{2n-i}: 1 \le i \le \frac{t}{2}-1\right\}$ is a right transversal of $H$ in $\langle x\rangle$. Since $x^{\frac{t(q-1)}{2}}\in H$, we have $x^n = x^{\frac{t(q-1)}{2}}x^{\frac{t}{2}} \in Hx^{\frac{t}{2}}$ and hence $A := \left\{e, x^n, x^i, x^{2n-i}: 1 \le i \le \frac{t}{2}-1\right\}$ is a right transversal of $H$ in $\langle x\rangle$. Moreover, $A$ is inverse-closed as $x^{2n} = e$. Let $B := \{x^i y: 0 \le i \le t-1\}$. Then $\cup_{b \in B} Hb = G\setminus\langle x\rangle$ and $Hb\ne Hb'$ for distinct $b, b' \in B$. Since $q \ge 3$ is odd, we have $x^{\frac{(q-1)t}{2}}\in H$. Hence $x^{n+j}y = x^{\frac{(q-1)t}{2}} (x^{\frac{t}{2}+j}y) \in H(x^{\frac{t}{2}+j}y)$ for $0\le j \le \frac{t}{2}-1$. It follows that $C := \left\{x^i y, x^{n+i} y: 0 \le i \le \frac{t}{2}-1\right\}$ satisfies $\cup_{c \in C} Hc = G\setminus\langle x\rangle$ and $Hc \ne Hc'$ for distinct $c, c' \in C$. Moreover, $C$ is inverse-closed as $(x^ry)^{-1}=x^{n+r}y$ for any integer $r$. Therefore, $A \cup C$ is an inverse-closed right transversal of $H$ in $G$. Setting $S := (A \cup C) \setminus\{e\}$, we obtain from Lemma~\ref{mainlem0} that $H$ is a perfect code in $\Cay(G, S)$. Note that $S$ is equal to the subset defined in \eqref{eq:s1}.
\qed

We can also prove the sufficiency of Lemma \ref{g1} using \cite[Lemma 2.10(a)]{HXZ18}. As usual, for a group $G$, let $\mathbb{Z}[G]$ be the group ring of $G$ over $\mathbb{Z}$. For a subset $A$ of $G$, write
$$
\overline{A}=\sum_{g\in G}\mu_A(g)g\in \mathbb{Z}[G],
$$
where
\begin{equation*}
\mu_A(g)=
\left\{
       \begin{array}{ll}
         1, & \hbox{$g\in A$;} \\
         0, & \hbox{$g\notin A$.}
       \end{array}
      \right.
\end{equation*}

\begin{lem}{\rm (\cite[Lemma 2.10(a)]{HXZ18})}
\label{basichuan}
Let $G$ be a group and ${\rm Cay}(G,S)$ a Cayley graph of $G$.
Let $C$ be a subset of $G$. Then $C$ is a perfect code in ${\rm Cay}(G,S)$ if and only if $\overline{S\cup\{e\}}\cdot\overline{C}=\overline{G}$.
\end{lem}

It is straightforward to verify that, for any positive divisor $t$ of $2n$ such that $\frac{2n}{t}$ is odd, we have $\overline{S\cup\{e\}}\cdot \overline{\langle x^t\rangle}=\overline{Q_{4n}}$, where $S$ is as defined in \eqref{eq:s1}. Hence, by Lemma \ref{basichuan}, $\langle x^t\rangle$ is a perfect code in $\Cay(Q_{4n}, S)$, proving the sufficiency of Lemma \ref{g1}. Using the same method, we can also prove the following lemma.

\begin{lem}\label{g2}
Let $n \ge 2$ be an integer, $t \ge 3$ an integer dividing $2n$, and $s$ an integer with $0\le s \le t-1$. Then the proper subgroup $\langle x^t, x^sy\rangle$ of $Q_{4n}$ is a perfect code
of $Q_{4n}$ if and only if $t$ is odd. Moreover, for any odd divisor $t \ge 3$ of $2n$ and any integer $s$ with $0\le s \le t-1$, $\Cay(Q_{4n}, S)$ with $S$ as given in \eqref{eq:s2} admits $\langle x^t, x^sy\rangle$ as a perfect code.
\end{lem}

\proof
Denote $G=Q_{4n}$ and $H=\langle x^t, x^sy\rangle$. Then $|H|$ is even and $|G:H|=t$. First, suppose that $H$ is a perfect code of $G$. That is, $H$ is a perfect code of some Cayley graph ${\rm Cay}(G,S)$ of $G$. Then, by Lemma~\ref{mainlem0}, $S\cup\{e\}$ is a left transversal of $H$ in $G$. Hence $|S|=t-1$. Suppose to the contrary that $t$ is even. Then $|S|$ is odd. Since $S^{-1}=S$, it follows that $x^n\in S$. On the other hand, as $|H|$ is even, we have $x^nH=H$, which contradicts the fact that $S\cup\{e\}$ is a left transversal of $H$ in $G$. Hence $t$ must be odd.

Conversely, suppose that $t \ge 3$ is odd. Note that
$$
\overline{H}=(e+x^t+x^{2t}+\cdots+x^{2n-t})(e+x^sy).
$$
Set $m = \frac{t-1}{2}$. Let $S=\{x,x^{-1},x^2,x^{-2},\ldots,x^m,x^{-m}\}$ be the subset of $G$ as given in \eqref{eq:s2}. Observe that $|S|=2m=t-1$.
It follows that
\begin{eqnarray*}
\overline{S\cup\{e\}}\cdot\overline{H}&=&(e+x+x^2+\cdots+x^{m}+
x^{-1}+x^{-2}+\cdots+x^{-m})
\\
& &(e+x^t+x^{2t}+\cdots+x^{2n-t})(e+x^sy)\\
&=&(e+x+x^2+\cdots+x^{2n-1})(e+x^sy)\\
&=&(e+x+x^2+\cdots+x^{2n-1})+(e+x+x^2+\cdots+x^{2n-1})x^sy\\
&=&\overline{G}.
\end{eqnarray*}
Therefore, by Lemma~\ref{basichuan}, $H$ is a perfect code of ${\rm Cay}(G,S)$.
\qed

\noindent {\em Proof of Theorem~{\rm\ref{newmainth2}}.}
Theorem~\ref{newmainth2} follows from Lemmas~\ref{g1} and \ref{g2} immediately.
\qed

We conclude the paper by the following example to illustrate Theorem \ref{newmainth2}.

\begin{ex}
\label{Q24}
Let $G=Q_{24}$. By Theorem~\ref{newmainth2}, we know that $G$, $\{e\}$, $\langle x^{4}\rangle$ and $\langle x^3,x^s y\rangle$ ($0\le s \le 2$) are all subgroup perfect codes of $G$. More explicitly, ${\rm Cay}(G, S_1)$ with $S_1=\{x,x^6,x^{11},y,xy,x^6y,x^7y\}$ admits $\langle x^{4}\rangle$ as a perfect code, and ${\rm Cay}(G, S_2)$ with $S_2=\{x,x^{11}\}$ admits $\langle x^3,x^s y\rangle$ as a perfect code for each $0\le s \le 2$. We draw ${\rm Cay}(G, S_1)$ in Figure~\ref{Q24-1}, where for each $i$ with $0 \le i \le 5$, $x^i$ is joined to $x^{i+6}$ by an edge. From this drawing one can easily see that $\langle x^{4}\rangle = \{e, x^4, x^8\}$ is a perfect code of ${\rm Cay}(G, S_1)$; that is, the vertices $e$, $x^4$ and $x^8$ are pairwise non-adjacent and every other vertex is adjacent to exactly one of these three vertices. One can also see that ${\rm Cay}(G, S_2)$ is disconnected with two connected components, namely the $12$-cycles $(e, x, x^2, \ldots, x^{11}, e)$ and $(y, xy, x^2 y, \ldots, x^{11}y, y)$.

\smallskip
\begin{figure}[hptb]
  \centering
  \includegraphics[width=10cm]{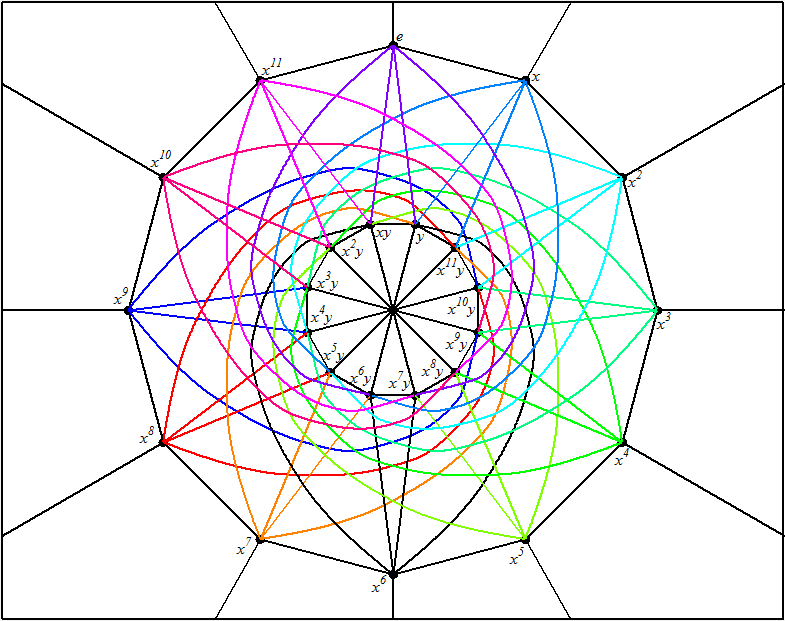}\\
  \caption{The graph ${\rm Cay}(G, S_1)$ in Example \ref{Q24}.}
  \label{Q24-1}
\end{figure}
\end{ex}

\bigskip
\noindent \textbf{Acknowledgements}~~We are grateful to the anonymous referees for their careful reading and helpful comments, and to Binzhou Xia for informing us an error in an early version of this paper. Ma was supported by the National Natural Science Foundation of China (Grant No. 11801441)
and the Natural Science Basic Research Program of Shaanxi (Program No. 2020JQ-761).
Wang was supported by the National Natural Science Foundation of China (Grant No. 11671043). Zhou was supported by the National Natural Science Foundation of China (Grant No. 61771019) and the Research Grant Support Scheme of The University of Melbourne.

\end{document}